\documentclass[12pt]{article}
\usepackage{latexsym}
\usepackage{amsmath,amsfonts,amssymb,mathrsfs,pslatex,graphicx}
\usepackage{amsthm}
\usepackage{verbatim}% for comment
\usepackage{pdfsync} 

\title{Multistable processes and localisability}
\author{Kenneth Falconer  and Lining Liu\\
\small{{\it Mathematical Institute,  
University of St~Andrews, North Haugh, St~Andrews,}} \\
\small{{\it Fife, KY16~9SS, Scotland }}} 
\def\bbbr{{\rm I\!R}}

\newcommand\ep{\stackrel{p}{=}}
\newcommand\topp{\stackrel{{\rm p}}{\rightarrow}}
\newcommand\tod{\stackrel{{\rm d}}{\rightarrow}}
\newcommand\fdd{\stackrel{{\rm fdd}}{\rightarrow}}

 \newtheorem{theo}{Theorem}[section]
 \newtheorem{cor}[theo]{Corollary}
 \newtheorem{lem}[theo]{Lemma}
 \newtheorem{prop}[theo]{Proposition}

 \newtheorem{exa}[theo]{Example}

\begin{document}
\maketitle

\begin{abstract} 
We use characteristic functions to construct $\alpha (x)$-multistable measures and integrals, where the measures behave locally like $\alpha$-stable measures, but with the stability index $\alpha (x)$ varying with time $x$. This enables us to construct  $\alpha (x)$-multistable processes on $\bbbr$, that is processes whose scaling limit at time $x$ is an $\alpha (x)$-stable process. We present several examples of such multistable processes and examine their localisability. 
\end{abstract}

 \section{Introduction}
 \setcounter{equation}{0}
\setcounter{theo}{0}
There are several ways of constructing $\alpha$-stable processes, that is stochastic processes  such that the finite dimensional distributions of the process at any finite set of $m$ times is an $m$-dimensional $\alpha$-stable vector, see \cite{Bk_Sam} for an full discussion. In this paper we construct $\alpha (x)$-multistable processes, that is processes which look locally like $\alpha (x)$-stable processes close to time $x$ in the sense that the local scaling limits are $\alpha (x)$-stable processes, but where the stability index $\alpha (x)$ varies with time $x$.

A number of constructions for multistable processes have been given recently, generalising the constructions of stable processes. One approach is based on Poisson point process \cite{FL}, and another is based on sums of random series  \cite{LL}. Here we use characteristic functions to construct multistable integrals and measures. We show that these multistable measures are locally like $\alpha$-stable measures and  may be approximated by sums of many independent $\alpha$-stable measures defined on short intervals with differing $\alpha$. We then use multistable integrals to define multistable processes and  give sufficient conditions for processes to be localisable or strongly localisable, that is to have a local scaling limit. We give a range of examples of such multistable processes.

 \section{Definition of $\alpha(x)$-multistable measure and integral}\label{ch3_sec2}
\setcounter{equation}{0}
\setcounter{theo}{0}
 Throughout this paper, for given  $0<a\leq b <\infty$, the function  $\alpha:\mathbb{R} \to [a,b]$ will be a Lebesgue measurable function that will play the r\^{o}le of a varying stability index. 
 We will work with various linear spaces of measurable functions on $\mathbb{R}$. For $0<p < \infty$ let
 $$ \mathcal{F}_{p} = \{ f: f \textrm{ is measurable with }  \|f\|_p <\infty\} 
  \textrm{ where } \|f\|_p= \bigg(\int|f(x)|^{p}dx\bigg)^{1/p};$$
  thus $ \| \cdot \|_p$ is a quasinorm (i.e. there is a weak triangle inequality $ \| f+g \|_p \leq k ( \| f \|_p+ \| g \|_p)$ for some $k>0$) which becomes a norm if $1 \leq p < \infty$.
It is convenient to  write
 $$|f(x)|^{a,b}=\max \big\{ |f(x)|^{a},
 |f(x)|^{b} \big \},$$
 and to define the space of the functions
 \begin{equation*}\label{Fab}
 \mathcal{F}_{a,b} =   \mathcal{F}_{a}\cap \mathcal{F}_{b}   =   \{ f: f \textrm{ is measurable with }
 \int|f(x)|^{a,b}dx<\infty \}.
\end{equation*}
We also define variable exponent Lebesgue spaces (special cases of Orlicz spaces, see for example \cite{ELN}) by
 $$ \mathcal{F}_{\alpha} = \{ f: f \textrm{ is measurable with }  \|f\|_\alpha <\infty\} 
  \textrm{ where } \|f\|_\alpha= \Big\{\lambda >0 :   \int\Big|\frac{f(x)}{\lambda}\Big|^{\alpha(x)}dx  = 1\Big\}.$$
  Then $ \|\cdot \|_\alpha$ is a quasinorm that reduces to $ \|\cdot \|_p$ if $\alpha(x)= p$ is constant, and is a norm if $1\leq a \leq \alpha (x) \leq b$ for all $x$.
  
  Note that with $a \leq \alpha(x) \leq b$ we have  $\mathcal{F}_{\alpha} \subseteq   \mathcal{F}_{a,b}$ with
  $$ \|f\|_\alpha \leq c_{a,b}\max \big\{ \|f(x)\|_{a},\|f(x)\|_{b} \big \},$$
  where $c_{a,b}$ depends only on $a$ and $b$.

We define the multistable stochastic integral $I(f)$ of a function $f \in \mathcal{F}_{\alpha}$ by specifying the
finite-dimensional distributions of $I$ as a stochastic process on the space of functions $\mathcal{F}_{a,b}$ and then using the Kolmogorov Extension
Theorem to show that the process is well-defined. 

Given $f_{1},f_{2},...,f_{d} \in \mathcal{F}_{\alpha}$, the following proposition shows that we can define a probability distribution on
 the vector $(I(f_{1}),I(f_{2}),...,I(f_{d}))\in \mathbb{R}^{d}$ by the characteristic function $\phi_{f_{1},...f_{d}}$ given by (\ref{def_pm_1}). The essential point is that $\alpha(x)$ may vary with $x$.
 
 \begin{lem}\label{Ch3_def_pms}
 Let $d \in \mathbb{N}$ and $f_{1}, f_{2},...,f_{d} \in \mathcal{F}_{\alpha}$, where $0<a \leq \alpha(x) \leq 2$ for all $x \in\mathbb{R}$ . Then 
 \begin{eqnarray}
 \nonumber \phi_{f_{1},...f_{d}}(\theta_{1},...,\theta_{d})&=&\mathbb{E} \Big( \exp i\sum_{j=1}^{d} \theta_{j} I(f_{j}) \Big)\\
 \label{def_pm_1} &=&\exp  \Big\{-\int
 \big| \sum_{j=1}^{d} \theta_{j} f_{j}(x)  \big|^{\alpha(x)}dx  \Big\}
 \end{eqnarray}
 for $(\theta_1, \theta_2, \dots, \theta_d) \in \mathbb{R}^{d}$, is the characteristic function of a probability distribution on the random vector $(I(f_{1}),I(f_{2}),...,I(f_{d}))$.
 \end{lem}
 
\begin{proof}
First, assume that $\alpha(x)$ is given by the simple function
\begin{equation}\label{prop-simple-1}
\alpha(x)= \sum_{k=1}^{m} \alpha_{k} \mathbf{1}_{A_k}(x),
\end{equation}
where $ 0 < \alpha_k \leq 2$ and $A_k$ are disjoint Lebesgue measurable sets with $\cup_{k=1}^{m}A_k= \mathbb{R}$.

For $(\theta_1, \dots, \theta_d) \in  \mathbb{R}^d$
\begin{eqnarray}
\exp \Big\{ -\int \big|\sum_{j=1}^{d} \theta_j f_j(x) \big|^{\alpha(x)}dx \Big\}&=&\exp \Big\{-\sum_{k=1}^{m} \int \big|\sum_{j=1}^{d} \theta_j f_j(x)\mathbf{1}_{A_k}(x) \big|^{\alpha_k}dx \Big\}\nonumber\\
 &=& \prod_{k=1}^{m} \exp \Big\{-\int  \big|\sum_{j=1}^{d} \theta_j f_j(x)\mathbf{1}_{A_k}(x) \big|^{\alpha_k}dx \Big\}.
\label{product-alpha-stable}
\end{eqnarray}
Now, $\exp\big\{-\int \big|\sum_{j=1}^{d} \theta_j f_j(x)\mathbf{1}_{A_k}(x)\big|^{\alpha_k}dx\big\}$ is the characteristic function of the $\alpha_k$-stable random vector $(I(f_1\mathbf{1}_{A_k}), \dots, I(f_d\mathbf{1}_{A_k}))$, see \cite{Bk_Sam}. Hence (\ref{product-alpha-stable}) is the product of the characteristic functions of $m$ $\alpha_k$-stable random vectors and so is the characteristic function of a $d$-dimensional random vector given by the independent sum of $\alpha_k$-stable random vectors. Hence (\ref{def_pm_1}) is a valid characteristic function of a random vector $(I(f_1), \dots, I(f_n))$ in the case when $\alpha(x)$ is a simple function (\ref{prop-simple-1}).

Now let $0 <a \leq  \alpha(x) \leq 2$ be measurable. Given   $f_1, \dots f_d \in \mathcal{F}_{\alpha}$  write $A=\{x: \sum_{j=1}^{d}|f_j(x)| \leq 1\} $.
Take a sequence of simple functions $\{\alpha_p(x)\}_{p=1}^{\infty}$ with $0 < \alpha_p(x) \leq 2$ such that $\alpha_p(x) \to \alpha(x)$ pointwise almost everywhere; we may assume that $\alpha_p(x) \geq \alpha(x)$ if $x \in A$ and   $\alpha_p(x) \leq \alpha(x)$ if $x \notin A$, for all $x$ and $p$.
Then
\begin{eqnarray*}
\big|\sum_{j=1}^{d} \theta_j f_j(x)\big|^{\alpha_p(x)}& \leq&\max_j \max\{|\theta_j|^a,|\theta_j|^2\} \big(\sum_{j=1}^{d} | f_j(x)|\big)^{\alpha_p(x)}\\
&\leq&\max_j \max\{|\theta_j|^a,|\theta_j|^2\} \big(\sum_{j=1}^{d} | f_j(x)|\big)^{\alpha(x)},
\end{eqnarray*}
an expression that is integrable since
$f_1, \dots f_d \in \mathcal{F}_{\alpha}$. By the dominated convergence theorem,
\begin{equation}\label{inter-conver}
\int \big|\sum_{j=1}^{d} \theta_j f_j(x) \big|^{\alpha_p(x)} dx \to \int  \big|\sum_{j=1}^{d} \theta_j f_j(x) \big|^{\alpha(x)}dx,
\end{equation}
and so 
\begin{equation}\label{exponent-conver}
\exp\Big\{\int  \big|\sum_{j=1}^{d} \theta_j f_j(x)\big|^{\alpha_p(x)} dx \Big\} \to \exp\Big\{\int  \big|\sum_{j=1}^{d} \theta_j f_j(x)\big|^{\alpha(x)}dx \Big\},
\end{equation}
as $p \to \infty$, for all $\theta_1, \dots ,\theta_d \in \mathbb{R}$. 

For $f_1, \dots, f_d \in \mathcal{F}_{\alpha}$, 
$$
\int \big|\sum_{j=1}^{d} \theta_j f_j(x) \big|^{\alpha(x)}dx \leq\max_j \max\{|\theta_j|^a,|\theta_j|^2\} \sum_{j=1}^{d} | f_j(x)|^{\alpha(x)}
\to 0
$$
as $\max_j\{|\theta_j |\} \to 0$. Thus (\ref{def_pm_1}) is continuous at 0. Moreover from (\ref{product-alpha-stable}) $\exp  \big\{-\int  \big|\sum_{j=1}^{d} \theta_j f_j(x) \big|^{\alpha_p(x)}dx \big\}$ is a valid characteristic function of a $d$-dimensional random vector for all $p$. Applying L\'{e}vy's continuity theorem  to (\ref{exponent-conver}), there is a probability distribution on the random vector $(I(f_{1}),I(f_{2}),...,I(f_{d}))$, with characteristic function given by (\ref{def_pm_1}).
\end{proof}

As with $\alpha$-stable  integrals, see \cite{Bk_Sam}, Kolmogorov's extension theorem allows us to define $\alpha(x)$-stable integrals consistently on $ \mathcal{F}_{\alpha}$.

 \begin{theo}\label{Ch3_prop_pms}
Let $0<a\leq \alpha(x)\leq 2$.  There exists a stochastic process $\{I(f), f\in
 \mathcal{F}_{\alpha}\}$ with finite-dimensional distributions 
 given by ($\ref{def_pm_1}$), that is with $\phi_{I(f_1),\dots,I(f_d)}=\phi_{f_1,\dots,f_d}$ for all $f_1, \dots, f_d \in \mathcal{F}_{\alpha}$.
 \end{theo}
 \begin{proof}
 For $f_1,\dots,f_d \in \mathcal{F}_{\alpha}$ it follows from (\ref{def_pm_1}) that, for any permutation $(\pi(1),\pi(2),...,\pi(d))$ of
 $(1,2,...,d)$, we have
$$
 \phi_{f_{\pi(1),...,\pi(d)}}(\theta_{\pi(1)},...,\theta_{\pi(d)})=
 \phi_{f_{1},...,f_{d}}(\theta_{1},...,\theta_{d}),
$$
 and also that, for any $n\leq d$,
$$ \phi_{f_{1},...,f_{n}}(\theta_{1},...,\theta_{n})=
\phi_{f_{1},...,f_{n},...,f_{d}}(\theta_{1},...,\theta_{n},0,...,0).
 $$
 Thus the probability distributions given by (\ref{def_pm_1}) satisfy the consistency conditions for  Kolmogorov's Extension Theorem, so, applying this theorem to the space of functions $\mathcal{F}_{\alpha}$, there is a stochastic process on $\mathcal{F}_{\alpha}$ which we denote by $\{I(f),f\in\mathcal{F}_{\alpha}\}$, whose
 finite-dimensional
 distributions are given by the characteristic functions (\ref{def_pm_1}).
 \end{proof}

 We call $I(f)$ the \textit{$\alpha(x)$-multistable integral} of $f$.
 By applying  (\ref{def_pm_1}) with functions $(a_1f_1+a_2 f_2)$, $f_1$, $f_2$ and variables $\theta$, $-a_1 \theta$, $-a_2 \theta$ it follows that the multistable integral is linear, that is if $f_{1},f_{2}\in \mathcal{F}_{\alpha}$ and $a_{1}, a_{2} \in \mathbb{R}$, then
 \begin{equation}\label{prop_line_1}
 I(a_{1}f_{1}+a_{2}f_{2})= a_{1}I(f_{1})+a_{2}I(f_{2}) \textrm{ } \textrm{
 }
 \textrm{ } \textrm{ }\textrm{ a.s. }
 \end{equation}
 
% With $(\Omega,F,\mathbb{P})$ the underlying
% probability space, we write $\mathcal{L}^{0}(\Omega)$ for the set of all real
% random variables defined on $\Omega$. 
 Let $L$ be Lebesgue measure on $\mathbb{R}$, let $\mathcal{E}$ be the Lebesgue measurable subsets of $\mathbb{R}$ and let
 $\mathcal{E}_{0}=\{ A\in \mathcal{E}:L(A)<\infty \}$
 be the sets of finite Lebesgue
 measure. Let $\alpha: \mathbb{R} \to [a,b]$ be measurable where $0<a\leq b\leq2$. Analogously to \cite[Section 3.3]{Bk_Sam} for $\alpha$-stable measures, we define the $\alpha(x)$-{\em multistable random measure} $M$ by
\begin{equation}\label{def_meas_1}
 M(A)=I(\mathbf{1}_{A})
 \end{equation} 
for $A \in  \mathcal{E}_{0}$, where $\mathbf{1}_{A}$ is the indicator function of the set $A$; thus $M(A)$ is a random variable for each $A \in  \mathcal{E}_{0}$.

  It is natural to write
 \begin{equation}
 \int f(x)dM(x): =I(f), \textrm{ } \textrm{ } \textrm{ } \textrm{
 } f\in \mathcal{F}_{\alpha}, \label{defint}
 \end{equation}
 since there are many analogues to usual integration with respect to a measure. Clearly, linearity of this integral is a restatement of  (\ref{prop_line_1}), and
 \begin{equation*}
 \int \mathbf{1}_{A}(x)dM(x)=M(A).
 \end{equation*}
   With this notation the characteristic function (\ref{def_pm_1}) may be written
 \begin{equation} \label{def_meas_33}
 \mathbb{E}\Big( \exp i\big\{ \sum_{j=1}^{d} \theta_{j} \int
 f_{j}(x)dM(x) \big\}\Big)=\exp \Big\{-\int \big|
 \sum_{j=1}^{d}\theta_{j}f_{j}(x) \big|^{\alpha(x)}dx\Big\}
 \end{equation}
 for $f_{j}\in \mathcal{F}_{\alpha}$. 
 For the random measures, taking $f_{j}=\mathbf{1}_{A_{j}}$ with
$A_{j}\in \mathcal{E}_{0}$,
\begin{equation}\label{def_meas_4}
\mathbb{E}\Big(\exp i\big\{\sum_{j=1}^{d}\theta_{j}M(A_{j})
 \big\} \Big)=\exp \Big\{-\int \big| \sum_{j=1}^{d}\theta_{j}\mathbf{1}_{A_{j}}(x)
 \big|^{\alpha(x)}dx \Big\}.
 \end{equation}

We may estimate the moments of an $\alpha(x)$-multistable integral in terms of the norm $\|\cdot\|_\alpha$.

\begin{prop}\label{moments}
Let $0 < a \leq \alpha(x) \leq b \leq 2$ and let $g \in \mathcal{F}_{\alpha}$. Then there is a number $c_1$ depending only on $a$ and $b$ such that for all $\lambda >0$
\begin{equation}
\mathbb{P}\bigg(\Big|\int g(x) dM_\alpha(x)\Big| \geq \lambda\bigg)
\leq c_1 \int\bigg| \frac{g(x)}{\lambda}\bigg|^{\alpha(x)}dx. \label{probbig}
\end{equation}
Moreover, if $0<p < \inf_{x \in \mathbb{R}} \alpha(x)$ there is a number $c_2$ depending only on $p$ and $b$ such that 
\begin{equation}
\mathbb{E}\bigg(\Big|\int g(x) dM_\alpha(x)\Big|^p\bigg) \leq c_2\|g\|_\alpha^p.  \label{moment}
\end{equation}
\end{prop}

\begin{proof}
A simple calculation using distribution functions (see 
\cite[p.47]{Bil}) gives
\begin{eqnarray*}
\mathbb{P}\bigg(\Big|\int g(x) dM_\alpha(x)\Big| \geq \lambda\bigg)& \leq &
    \frac{\lambda}{2}\int_{-2/\lambda}^{2/\lambda}\bigg(1-\mathbb{E}\Big(\exp\big(i\theta    \int g(x) dM_{\alpha(x)}\big)\Big)\bigg)d\theta 
    \\
    &= & 
    \frac{\lambda}{2}\int_{-2/\lambda}^{2/\lambda}
    \Big(1-\exp\big(- \int |\theta g(x)|^{\alpha(x)}dx \big)\Big) d\theta\\
     &\leq & 
    \frac{\lambda}{2}\int_{-2/\lambda}^{2/\lambda}
    \Big( \int |\theta|^{\alpha(x)} |g(x)|^{\alpha(x)}dx\Big)d\theta\\
     &\leq & 
c_1 \int\Big| \frac{g(x)}{\lambda}\Big|^{\alpha(x)}dx. 
\end{eqnarray*}

Assuming as we may that $c_1 \geq 1$ and writing $\lambda_0 =\|g\|_{\alpha} >0$ for the number such that 
$ \int \lambda_0^{-\alpha(x)}|g(x)|^{\alpha(x)}dx=1$, we have
\begin{eqnarray*}
\mathbb{E}\bigg(\Big|\int g(x) dM_\alpha(x)\Big|^p\bigg) 
&=& p\int_{0}^{\infty}\lambda^{p-1}\mathbb{P}\bigg(\Big|\int g(x) dM_\alpha(x)\Big| \geq \lambda\bigg)d \lambda\\
&\leq &  c_1 p\int_{0}^{\infty}\lambda^{p-1} \min \Big\{1,\int \lambda^{-\alpha(x)}|g(x)|^{\alpha(x)}dx \Big\}d \lambda\\
 &\leq &    c_1 p \int_{0}^{\lambda_0}\lambda^{p-1}d\lambda 
    +   c_1 p  \int\int_{\lambda_0}^{\infty}\lambda^{p-1-\alpha(x)}|g(x)|^{\alpha(x)}
d\lambda dx 
\\
& \leq & 
c_3 \lambda_0^p + c_3  \lambda_0^p \int  \lambda_0^{-\alpha(x)}|g(x)|^{\alpha(x)}
dx \\
&=&
 c_2\|g\|_{\alpha}^{p}.
\end{eqnarray*}    
\end{proof}

Recall that a random measure $M$ on $\mathbb{R}$ is {\em independent scattered} if
$M(A_{1}), M(A_{2}),..., M(A_{d})$ are independent whenever  $A_{1}$,$A_{2}$,...,$A_{k} \in \mathcal{E}_{0}$ are pairwise disjoint, and is $\sigma${\em -additive} if whenever $A_{1}$,$A_{2}$,...$\in$
 $\mathcal{E}_{0}$ are disjoint and
 $\bigcup_{j=1}^{\infty}A_{j}\in \mathcal{E}_{0}$ then almost surely
 $$M \big(\bigcup_{j=1}^{\infty}A_{j}\big)={\sum_{j=1}^{\infty}}M(A_{j}),$$ taking an independent sum.

 \begin{theo}\label{ch3_theo_inscat}
 The $\alpha(x)$-multistable measure $M_{\alpha}$ is independent scattered and $\sigma$-additive.
 \end{theo}
 \begin{proof}
 This is a slight variant of  \cite[Section 3.3]{Bk_Sam}.
 Let $A_{1}$,$A_{2}$,...,$A_{k}$ $\in$ $\mathcal{E}_{0}$ be pairwise disjoint. Then using  (\ref{def_meas_4})
 \begin{equation*}
 \mathbb{E}\Big(\exp \big\{ i\sum_{j=1}^{d}\theta_{j}M_{\alpha}(A_{j})
 \big\} \Big) 
=\prod_{j=1}^{d}\exp \Big\{ -\int \big|
 \theta_{j}\mathbf{1}_{A_{j}}(x)
 \big|^{\alpha(x)}dx \Big\}\\
=\prod_{j=1}^{d}\mathbb{E}(\exp\{i\theta_{j}M_{\alpha}(A_{j})\}).
 \end{equation*}
 so $M_{\alpha}(A_{1}), M_{\alpha}(A_{2}),..., M_{\alpha}(A_{d})$ are independent, and $M_{\alpha}$ is
 independent scattered. 

If $A_{1}$,$A_{2}$,...,$A_{k}$ $\in$ $\mathcal{E}_{0}$ is a finite collection of disjoint sets, using (\ref{def_meas_1}) and (\ref{prop_line_1}),
 $$
 M_{\alpha}\Big(\bigcup_{j=1}^{k}A_{j}\Big) = I\big(\mathbf{1}_{\cup_{j=1}^{k}{A_{j}}}\big)
 = I\Big(\sum_{j=1}^{k}\mathbf{1}_{A_{j}}\Big)
 = \sum_{j=1}^{k}I(\mathbf{1}_{A_{j}})
 = \sum_{j=1}^{k}M_{\alpha}(A_{j}).
$$
 For a countable family of disjoint sets $A_{1},A_{2},... \in \mathcal{E}_{0}$ with
$B\equiv \bigcup_{j=1}^{\infty}A_{j}\in \mathcal{E}_{0}$, so that  $B=\bigcup_{j=1}^{k}A_{j} \cup \Big(\bigcup_{j=k+1}^{\infty}A_{j}\Big)$,
 it follows from above
 that
 \begin{equation}
M_{\alpha}(B)=M_{\alpha}\Big(\bigcup_{j=1}^{k}A_{j}\Big)+M_{\alpha}\Big(\bigcup_{j=k+1}^{\infty}A_{j}\Big)
= \sum_{j=1}^{k}M_{\alpha}(A_{j})+M_{\alpha}\Big(\bigcup_{j=k+1}^{\infty}A_{j}\Big). \label{thm_add_1} 
 \end{equation}
 Since
 $\lim_{k\to\infty}\mathcal{L}(\bigcup_{j=k+1}^{\infty}A_{j})=0$ and $\alpha(x) \in [a,b]$, for each $\theta \in \bbbr$
 $$
 \lim_{k\to\infty}\mathbb{E}\Big(\exp i\big\{ \theta
 M_{\alpha}(\bigcup_{j=k+1}^{\infty}A_{j}) \big\}
 \Big)=\lim_{k\to\infty}\exp\Big\{-\int|\theta\mathbf{1}_{\cup_{j=k+1}^{\infty}A_{j}}|^{\alpha(x)}\Big\} = 1,
$$
so $M_{\alpha}(\bigcup_{j=k+1}^{\infty}A_{j}) \tod 0$
 as $k\to \infty$ by L\'{e}vy's Continuity Theorem.

 By (\ref{thm_add_1}) we get $M_{\alpha}(B)-\sum_{j=1}^{k}M_{\alpha}(A_{j})\tod 0$ and
 so $M_{\alpha}(B)-\sum_{j=1}^{k}M_{\alpha}(A_{j})\topp 0$ as $k\to \infty$. Thus $\lim_{k \to \infty}\sum_{j=1}^{k}M_{\alpha}(A_{j})\ep M_{\alpha}(B)$, and, since the summands
 $M_{\alpha}(A_{j})$ are independent, this implies convergence almost surely, by a theorem of Kolmogorov, see \cite{Bk_LT}.
Thus $M_{\alpha}$
 is $\sigma$-additive.
 \end{proof}

Next we obtain conditions for convergence of a sequence of multistable measures with different multistable indexes.

\begin{prop}\label{ch3_thm_convergence}
 Let $\alpha_{n}(x)$, $\alpha(x)$ be Lebesgue measurable with
 $0<a\leq \alpha_{n}(x), \alpha(x) \leq b \leq 2$ for all $x \in
 \mathbb{R}$. Let $M_{\alpha_n}, M_{\alpha}$ be the associated $\alpha_{n}(x)$-multistable
 and $\alpha(x)$-multistable measures characterised by (\ref{def_meas_4}). Suppose $\alpha_{n}(x) \to \alpha(x)$
 for almost all $x \in \mathbb{R}$. Then
 $M_{\alpha_n} \fdd M_{\alpha}$ as $n\to\infty$, that is for
 all $m \in \mathbb{N}$ and $A_{1},A_{2},...,A_{m} \in \mathcal{E}_{0}$,\\
 $$(M_{\alpha_n}(A_{1}),
 M_{\alpha_n}(A_{2}),...,M_{\alpha_n}(A_{m}))\tod (M_{\alpha}(A_{1}),
 M_{\alpha}(A_{2}),...,M_{\alpha}(A_{m})).$$ 
 \end{prop}
 \begin{proof}
 Let $A_{1}$, $A_{2}$,...,$A_{m}$ $\in \mathcal{E}_{0}$.
  Then for all $n$ and all $x\in \mathbb{R}$
 $$\big|\sum_{j=1}^{m} \theta_{j}\mathbf{1}_{A_{j}}(x)\big|^{\alpha_{n}(x)} \leq c\mathbf{1}_{A}(x)$$
where $A=\bigcup_{j=1}^{m}A_{j} \in \mathcal{E}_{0}$ and
  $c=\max\big\{ \big(\sum_{j=1}^{m}|\theta_{j}|\big)^{a},  \big(\sum_{j=1}^{m}|\theta_{j}|\big)^{b} \big\}.$
Since $\int \mathbf{1}_{A}(x)dx < \infty$, 
the dominated convergence theorem implies that
 $$\lim_{n\to \infty}\exp\Big(-\int\big|\sum_{j=1}^{m}\theta_{j}\mathbf{1}_{A_{j}}(x)\big|^{\alpha_{n}(x)}dx\Big)=\exp\big(-\int\Big|\sum_{j=1}^{m}\theta_{j}\mathbf{1}_{A_{j}}(x)\big|^{\alpha(x)}dx\Big),$$ 
so from (\ref{def_meas_4}), 
 $$ \mathbb{E}\Big(\exp i \big\{\sum_{j=1}^{m}\theta_{j}M_{\alpha_n}(A_{j})
 \big\}
 \Big) \to  \mathbb{E}\Big(\exp i \big\{\sum_{j=1}^{m}\theta_{j}M_{\alpha}(A_{j})
 \big\}\Big)$$
 as $n\to \infty$.
  By L\'{e}vy's continuity theorem  $M_{\alpha_n} \fdd M_{\alpha}$.
 \end{proof}

 To get a feel for $\alpha(x)$-mutistable measures, we show that, for a continuous $\alpha(x)$, the $\alpha(x)$-multistable measure $M$ may be approximated by random measures that are the sum of many independent 
 $\alpha$-stable measures defined on short intervals.
 
 Assume  that $\alpha : \mathbb{R} \to [a,b] \subset (0,2]$ is continuous and let $M_{\alpha}$ be the $\alpha(x)$-multistable measure on the sets $\mathcal{E}_{0}$. We now use the same procedure but with piecewise constant functions  $\alpha_{n}(x)$ to obtain approximating measures $M_{\alpha_n}$. 
 
 For each $n$ let $\alpha_{n}: \mathbb{R} \to [a,b] \subset (0,2)$ be given by
 $$\alpha_{n}(x)=\alpha(r2^{-n}) \mbox{ if } x \in [r2^{-n}, (r+1)2^{-n}) \mbox{ for } 
 r \in \mathbb{Z}$$
and let  $M_{\alpha_n}$ be the $\alpha_n$-multistable measure obtained from $\alpha_n(x)$ as above, so in particular  $M_{\alpha_n}$ has finite-dimensional distributions given by the characteristic function
 \begin{equation*}\label{thm_chf_1}
 \mathbb{E} \Big(\exp i  \big\{\sum_{j=1}^{d}\theta_{j}M_{\alpha_n}(A_{j})
 \big\}  \Big)
 =\exp \Big\{-\int \big| \sum_{j=1}^{d} \theta_{j}
 \mathbf{1}_{A_{j}}(x) \big|^{\alpha_{n}(x)}dx \Big\}.
 \end{equation*}
It follows from Theorem \ref{ch3_theo_inscat} that each $M_{\alpha_n}$ is independent scattered and $\sigma$-additive.

\begin{theo}
Let $0<a\leq b\leq 2$ and $\alpha: \mathbb{R} \to [a,b]$ be continuous.  Let $M_{n,r}$ denote the restriction of $\alpha(r2^{-n})$-stable measure to the interval $[r2^{-n},(r+1)2^{-n}))$, that is
\begin{equation*}\label{def_vmn_1}
 M_{n,r}(A)= M_{\alpha(r2^{-n})}  (A\cap [r2^{-n},(r+1)2^{-n}))=M_{\alpha_n(x)}  (A\cap [r2^{-n},(r+1)2^{-n})),
 \end{equation*}
where $M_{\alpha(r2^{-n})}$ is $\alpha(r2^{-n})$-stable measure.
Then  $M_{\alpha_n}$ is a random measure given by  the independent sum of random measures
 \begin{equation*}\label{theo_vmn_1}
 M_{\alpha_n}(A)=\sum_{r\in\mathbb{Z}}M_{n,r}(A)
 \end{equation*}
 almost surely for $A \in \mathcal{E}_{0}$.
 Moreover $M_{\alpha_n} \fdd M_{\alpha}$ as $n\to \infty$.
\end{theo}

\begin{proof}
Since $M_{\alpha_n}$ is independent scattered, we have that for each $A \in \mathcal{E}$
$$M_{\alpha_n}(A\cap [r2^{-n},(r+1)2^{-n}))= M_{n,r}(A)$$
are independent for distinct $r$.

Let $A \in \mathcal{E}_{0}$. Since $M_{\alpha_n}$ is $\sigma$-additive,
\begin{eqnarray*}
M_n(A)&=&M_{\alpha_n}(A)\\
&=&M_{\alpha_n}\big(\bigcup_{r \in \mathbb{Z}}A\cap [r2^{-n},(r+1)2^{-n})\big)\\
&=&\sum_{r \in \mathbb{Z}}M_{\alpha_n}(A\cap [r2^{-n},(r+1)2^{-n}))\\
&=&\sum_{r \in \mathbb{Z}}M_{n,r}(A)
\end{eqnarray*}
where the summands are independent.

 For  each $n$ we have $\alpha_{n}(x)=\alpha(r2^{-n})$ for all $x \in
 [r2^{-n},(r+1)2^{-n})$. Since  $\alpha(x)$ is assumed continuous, we have
 $\lim_{n\to \infty}\alpha_{n}(x)=\alpha(x)$
 for all $x$. Thus by Theorem \ref{ch3_thm_convergence},  $M_{\alpha_n} \fdd M_{\alpha}$ as $n\to \infty$.
\end{proof}

One would expect an $\alpha(x)$-multistable
 measure to `look like' an $\alpha(u)$-stable measure in a small interval around $u$. We now make this idea precise.

 For $u \in \mathbb{R}$, $r>0$, let $T_{u,r}: \mathbb{R}\to \mathbb{R}$ be
 the scaling map, $T_{u,r}(x)=(x-u)/r.$
 This induces a mapping $T_{u,r}^{\#}$ on random integrals and measures, given by
 \begin{equation}
 \nonumber \int f(x)d(T_{u,r}^{\#}M_{\alpha})(x) = \int f\Big(\frac{x-u}{r}\Big)dM_{\alpha}(x)
\equiv I\Big(f\Big(\frac{.-u}{r}\Big)\Big). \label{def_map_1}
 \end{equation}
In particular 
$$(T_{u,r}^{\#}M_{\alpha})(A)=M_{\alpha}(T_{u,r}^{-1}(A))=I(\mathbf{1}_{T_{u,r}^{-1}(A)})$$
 for $A \in \mathcal{E}_{0}$
by (\ref{def_meas_1}).

 We show that scaling an $\alpha(x)$-multistable random measure about a point $u$ yields the $\alpha(u)$-stable measure $M_{\alpha(u)}$.

  \begin{theo}
 Let $\alpha: \mathbb{R}\to[a,b]\subseteq (0,2]$ be
 continuous with 
 \begin{equation}
 |\alpha(x+r)-\alpha(x)|=o(1/\log r)\label{contcond}
 \end{equation}
 uniformly on bounded intervals and let $u \in \mathbb{R}$. Then for all functions $f_1,\dots,f_d \in \mathcal{F}_{a,b}$ with compact support, the vectors
 \begin{eqnarray}
 \nonumber \Big(r^{-1/\alpha(u)}\int f_1(x)d(T_{u,r}^{\#}M_{\alpha})(x),\dots, r^{-1/\alpha(u)}\int f_d(x)d(T_{u,r}^{\#}M_{\alpha})(x)\Big)\hspace{2cm}\\
\label{theo-local-form-3} \tod \Big(\int f_1(x)dM_{\alpha(u)}(x),\dots, \int f_d(x)dM_{\alpha(u)}(x)\Big)
 \end{eqnarray}
 as $r\to 0$. In particular,
 \begin{equation}\label{theo-local-form-1}
r^{-1/\alpha(u)}\big((T_{u,r}^{\#}M_{\alpha})(A_1),\dots, (T_{u,r}^{\#}M_{\alpha})(A_d)\big)\tod (M_{\alpha(u)}(A_1),\dots, M_{\alpha(u)}(A_d))
\end{equation}
as $r\to 0$, for all bounded sets $A_1, \dots, A_d \in \mathcal{E}_{0}$. 
 \end{theo}
 \begin{proof}
  Let $f_{1},f_{2},...,f_{m} \in \mathcal{F}_{a,b}$ be functions with compact support, say
 in $[-z_{0},z_{0}]$. Let $\theta_{j}\in \mathbb{R}$,
 $j=1,2,...,m$, and consider the characteristic
 functions.
 \begin{eqnarray}
\nonumber &&\hspace{-2cm} \mathbb{E}\Big(\exp
 i\sum_{j=1}^{m}\theta_{j}r^{-1/\alpha(u)}\int
 f_{j}(x)d(T_{u,r}^{\#}M_{\alpha})(x)\Big) \\
 &=&\mathbb{E}\Big(\exp
 i\sum_{j=1}^{m}\theta_{j}r^{-1/\alpha(u)}\int f_{j}\Big(\frac{x-u}{r}\Big)dM_{\alpha}(x)\Big)\nonumber \\
\nonumber   &=&\exp\Big(-\int
 \big|\sum_{j=1}^{m}\theta_{j}r^{-1/\alpha(u)}f_{j}\Big(\frac{x-u}{r}\Big)\big|^{\alpha(x)}dx\Big) \\
 \nonumber  &=&\exp\Big(-\int \big|\sum_{j=1}^{m}\theta_{j}r^{-1/\alpha(u)}f_{j}(z)\big|^{\alpha(rz+u)}rdz\Big)\\
 &=&\exp\Big(-\int
 \big|\sum_{j=1}^{m}\theta_{j}f_{j}(z)\big|^{\alpha(rz+u)}r^{1-\alpha(rz+u)/\alpha(u)}dz\Big),\label{thm_local_1}
 \end{eqnarray}
 on writing $ (x-u)/r = z$.
 
From condition (\ref{contcond}) it is easy to see that
 $\lim_{r \to 0} r^{1-\alpha(rz+u)/\alpha(u)}=1$
 uniformly for $z \in [-z_{0},z_{0}]$, and also $\lim_{r\to
 0}\alpha(rz+u)=\alpha(u)$ uniformly for all $z \in [-z_{0},z_{0}]$ since $\alpha$ is continuous.
Noting that
there is a constant $c$ such that  for $r$ sufficiently small,
 $$\big|\sum_{j=1}^{m}\theta_{j}f_{j}(z)\big|^{\alpha(rz+u)}r^{1-\alpha(rz+u)/\alpha(u)}\leq
 c\sum_{j=1}^{m}|f_{j}(z)|^{a,b},$$
 for $z \in [-z_0,z_0]$ and $f_{j}\in \mathcal{F}_{a,b}$,  the
 dominated convergence theorem gives
 $$\lim_{r \to 0} \exp\Big(-\int
 \big|\sum_{j=1}^{m}\theta_{j}f_{j}(z)\big|^{\alpha(rz+u)}r^{1-\alpha(rz+u)/\alpha(u)}dz\Big) = \exp\big(-\int
\big|\sum_{j=1}^{m}\theta_{j}f_{j}(x)\big|^{\alpha(u)}dx\big),$$
so by (\ref{thm_local_1})
$$
\lim_{r\to 0}\mathbb{E}\Big(\exp
 i\sum_{j=1}^{m}\theta_{j}r^{-1/\alpha(u)}\int
 f_{j}(x)d(T_{u,r}^{\#}M_{\alpha})(x)\Big)
 =\mathbb{E}\Big(\exp
 i\sum_{j=1}^{m}\theta_{j}\int f_{j}(x)dM_{\alpha(u)}(x)\Big).
$$
L\'{e}vy's continuity theorem now implies 
 (\ref{theo-local-form-3}) and  (\ref{theo-local-form-1}).
 \end{proof}

\section{Multistable processes and localisability}
\setcounter{equation}{0}
\setcounter{theo}{0}
In this section we introduce processes defined by  multistable integrals, and in particular consider their local form, with the aim of  constructing processes with a prescribed local form. Thus, given $\alpha :  \mathbb{R} \to [a,b] \subset (0,2]$, we write
\begin{equation}\label{ydef}
Y(t)=\int f(t,x)dM_{\alpha}(x),
\end{equation}
for $t \in \mathbb{R}$ and $f \in \mathcal{F}_{a,b}$,
where the integrals are with respect to an  $\alpha(x)$-multistable measure $M_{\alpha}$ as in (\ref{defint}). By (\ref{def_meas_33}), for each $(t_1, t_2, \dots, t_d) \in \mathbb{R}^d$, the characteristic function of the random vector $(Y(t_1), Y(t_2), \dots, Y(t_d))$ is
\begin{eqnarray*}
\mathbb{E}\big(\exp i\sum_{j=1}^{d}\theta_j Y(t_j)\big)&=&\mathbb{E}\big(\exp i \sum_{j=1}^{d} \int f(t_j, x) dM_{\alpha}(x)\big)\\
&=&\exp\Big(-\int\big|\sum_{j=1}^{d} \theta_j f(t_j,x)\big|^{\alpha(x)}dx \Big).
\end{eqnarray*}
for all  $(\theta_1, \theta_2, \dots, \theta_d) \in \mathbb{R}^d$.

First we give conditions for $Y$ to have a continuous version.

\begin{prop}
Let    $\alpha: \mathbb{R} \to [a,b] \subseteq (1,2]$ be measurable and suppose $f(t,\cdot) \in {\cal F}_{\alpha}$ for all $t \in \mathbb{R}$. Let $Y$ be given by (\ref{ydef}). Suppose that there exists $1/a<\eta<1$, such that for each  bounded interval $I$ we can find $c>0$ such that
\begin{equation}
\|f(t,\cdot) -f(u,\cdot)\|_\alpha \leq c|t-u|^\eta \quad (t,u \in I). \label{ctsversion}
\end{equation}
Then  $Y$ has a continuous version satisfying an a.s. $\beta$-H\"{o}lder condition on each bounded interval for all $0< \beta <(\eta a -1)/a$.
In particular, (\ref{ctsversion}) holds if 
\begin{equation}
\int\big|
f(t,x)-f(u,x)\big|^{\alpha(x)}dx\leq c_1 |t-u|^{a\eta} \quad (t,u \in I),\label{ctsverequiv}
\end{equation}
a form that may be easier to check in practice.
\end{prop}

\begin{proof}
Take $p$ such that  $1/\eta <p <a$. By Proposition \ref{moments}
$$\mathbb{E}\big(|Y(t)-Y(u)|^p\big) 
= \mathbb{E}\bigg(\Big|\int(f(t,x)-f(u,x)) dM_\alpha (x)\Big|^p\bigg) 
\leq c_2\|f(t,\cdot) -f(u,\cdot)\|_\alpha^p \leq  c_2  c|t-u|^{\eta p}.$$
The conclusion follows from Kolmogorov's continuity theorem, see \cite[Theorem 25.2]{RW}. 
\end{proof}

Recall  that a stochastic process $Y$ is  localisable at a point if it has a
unique non-trivial scaling limit, formally $Y=\{Y(t):t\in \mathbb{R}\}$ is $h${\it -localisable} at $u$ with {\it local
form} or {\it tangent process} $Y_{u}'=\{Y_{u}'(t): t \in \mathbb{R}\}$ if
\begin{equation}
\frac{Y(u+rt)-Y(u)}{r^{h}} \fdd Y_{u}'(t) \label{converge}
\end{equation}
as $r \to 0$.
If $Y$ and  $Y_{u}'$ have versions in  $C(\bbbr)$ (the space of continuous
functions on $\bbbr$) and  convergence in (\ref{converge}) occurs in distribution with respect to the metric of uniform convergence on bounded intervals we say that $Y$ is {\it strongly localisable}. For the simplest example, a self-similar process with stationary 
increments $Y$ is 
localisable at all $u$ with local form $Y_{u}'=Y$ and is strongly localisable if it has a version in $C(\bbbr)$. In general there are considerable restrictions on the possible local forms, see \cite{Fal6}.

We call
a stochastic process $\{Y(t),\, t \in \mathbb{R} \}$ {\em multistable} if for almost all $u$, $Y$ is localisable at $u$ with $Y'_{u}$ an $\alpha$-stable process for some $\alpha=\alpha(u)$, where $0<\alpha(u)\leq 2$. Various constructions of multistable processes are given in  \cite{FL,FLL,LL}.

For a stochastic process $Y$, it is natural to ask under what
conditions $Y$ is localisable. The following theorem, which is a multistable analogue of \cite[Proposition 2.1]{FLL}, gives a sufficient
condition.

\begin{theo}\label{thm-main}
Let
\begin{equation}\label{msrd-def-si}
Y(t)= \int f(t,x)dM_{\alpha}(x),
\end{equation}
where $M_{\alpha}$ is an $\alpha(x)$-multistable measure for continuous $\alpha: \mathbb{R} \to [a,b] \subseteq (0,2]$. Assume that $f(t,.)\in
\mathcal{F}_{a,b}$ for all $t$ and
\begin{equation}\label{thm-main-condition-1}
\lim_{r\to 0}\int\Big|
\frac{f(u+rt,u+rz)-f(u,u+rz)}{r^{h-1/\alpha(u+rz)}}-h(t,z)\Big|^{a,b}dz=0
\end{equation}
for a jointly measurable function $h(t,z)$ with $h(t,.)\in
\mathcal{F}_{a,b}$ for all $t$. Then $Y$ is h-localisable at u with local form
\begin{equation}\label{thm-main-Yu}
Y'_{u}=\Big\{\int h(t,z) d M_{\alpha(u)}(z): t\in \mathbb{R}\Big\}
\end{equation}
where $M_{\alpha(u)}$ is $\alpha(u)$-stable measure.

If, in addition, there exists $\eta >1/a$ such that for each bounded interval $I$ we can find $c>0$ such that
\begin{equation}\label{stongcond}
\Big\|
\frac{f(u+rt,\cdot)-f(u+rv,\cdot)}{r^{h}}\Big\|_{\alpha}
\leq c|t-v|^\eta \quad (t,v \in I)
\end{equation}
for all sufficiently small $r>0$, then $Y$ is strongly localisable at $u$.
Condition (\ref{stongcond}) is implied by
\begin{equation}\label{equivstongcond}
\int\Big|
\frac{f(u+rt,u+rz)-f(u+rv,u+rz)}{r^{h-1/\alpha(u+rz)}}\Big|^{\alpha(u+rz)}dz\leq c_1 |t-v|^{a\eta}
\quad (t,v \in I)
\end{equation}
which can be more convenient to use in practice.
\end{theo}

To prove Theorem \ref{thm-main}, we need some convergence estimates.

\begin{lem}\label{lem3}
Let $0<a \leq b$. There is a constant $c$ that depends only on $a$ and $b$ such that, for all measurable $\alpha: \mathbb{R} \to [a,b]$ and $g,k\in \mathcal{F}_{a,b}$, 
\begin{align}\label{gkin}
\Big|\int &  |g(x)|^{\alpha(x) }dx  -\int |k(x)|^{\alpha(x)}dx \Big| \nonumber \\
\leq
& c \Big(\|g-k\|_a\|k\|_a^{\max\{0,a-1\}} +\|g-k\|_a^a
+\|g-k\|_b\|k\|_b^{\max\{0,b-1\}} +\|g-k\|_b^b \Big).
\end{align}
\end{lem}
\begin{proof}
If $0<a \leq \alpha(x) \leq b \leq 1$ for all $x \in \mathbb{R}$, then 
$$\Big|\int   |g(x)|^{\alpha(x) }dx  -\int |k(x)|^{\alpha(x)}dx \Big| 
\leq  \int |g(x)- k(x)|^{\alpha(x)}dx
\leq\|g-k\|_a^a + \|g-k\|_b^b.$$

On the other hand, if $1\leq a \leq \alpha(x) \leq b$ for all $x \in \mathbb{R}$, then by the mean value theorem there exists $0<\lambda(x) <1$ such that
\begin{eqnarray*}
\big| |g(x)|^{\alpha(x) }  - |k(x)|^{\alpha(x)} \big|
&  =  &  \alpha(x)  \big| |g(x)|- |k(x)|\big|\, \big|  |k(x)|+ \lambda(x)( |g(x)|-|k(x)|)\big|^{\alpha(x) - 1}\\
&\leq  & b  \big| |g(x)|- |k(x)|\big| \, \big|  |k(x)|+ |g(x)-k(x)|)\big|^{a-1}\\
&& + b  \big| |g(x)|- |k(x)|\big| \, \big|  |k(x)|+ |g(x)-k(x)|)\big|^{b-1}.
\end{eqnarray*} 
Integrating and using H\"{o}lder's inequality gives
%\begin{align*}
$$\Big|\int  |g(x)|^{\alpha(x) }dx  -\int |k(x)|^{\alpha(x)}dx \Big| \nonumber
\leq
 c \Big(\|g-k\|_a\big\| |k| +|g-k|\big\|_a^{a-1}
+\|g-k\|_b\big\| |k| +|g-k|\big\|_b^{b-1}\Big),$$
%\end{align*} 
which gives (\ref{gkin}) in this case.

In general, for $0 < a \leq \alpha(x) \leq b$, letting $A = \{x : a \leq \alpha(x) \leq 1\}$,  inequality (\ref{gkin}) holds for $g{\bf 1}_A$ and  $k{\bf 1}_A$ and also for $g{\bf 1}_{ \mathbb{R} \setminus A}$ and $k{\bf 1}_{ \mathbb{R} \setminus A}$,
 and combining these cases we get  (\ref{gkin}) for $g$ and $k$ for an appropriate $c$.
\end{proof}

We require the following Corollary.
\ref{thm-main}.

\begin{cor}\label{prop1}
Let $0<a \leq b$ and $g: \mathbb{R}^{+}\times \mathbb{R}\to
\mathbb{R}^*$ with $g(r,.)\in \mathcal{F}_{a,b}$ for all $r>0$. Let $k
\in \mathcal{F}_{a,b}$ and let $\beta:\mathbb{R}\to[a,b]$ be
continuous at $0$. If
\begin{equation}\label{prop1-condition-1}
\lim_{r\to 0} \int|g(r,z)-k(z)|^{a,b}dz= 0,
\end{equation}
then
\begin{equation}\label{prop1-condition-2}
\lim_{r\to 0}\int|g(r,z)|^{\beta(rz)}dz= \int|k(z)|^{\beta(0)}dz.
\end{equation}
\end{cor}
\begin{proof}
By (\ref{prop1-condition-1}) and Lemma \ref{lem3}
$$\lim_{r\to 0}\Big|\int   |g(r,z)|^{\beta(rz) }dz  -\int |k(z)|^{\beta(rz) }dz \Big| = 0.$$
Since $k \in \mathcal{F}_{a,b}$, the dominated convergence theorem gives
$$\lim_{r\to 0}\Big|\int |k(z)|^{\beta(rz) }dz  -\int |k(z)|^{\beta(0) }dz \Big| = 0,$$
and (\ref{prop1-condition-2}) follows on combining these two limits.
\end{proof}

We can now complete the proof of Theorem $\ref{thm-main}$.

\medskip

\noindent{\it Proof of Theorem \ref{thm-main}}
\quad Fix $u\in \mathbb{R}$. We consider the characteristic function of
the finite-dimensional distributions of $r^{-h}(Y(u+rt)-Y(u))$. Let
$\theta_{j}\in \mathbb{R}$ and $t_{j}\in \mathbb{R}$ for
$j=1,2,...,m$. Then, using (\ref{msrd-def-si}) and (\ref{def_meas_33}),
\begin{eqnarray}
\label{thm-main-proof-0}\mathbb{E} && \hspace{-1.0cm}\Big(\exp
i\sum_{j=1}^{m}\theta_{j}r^{-h}(Y(u+rt_{j})-Y(u))\Big)\\
\nonumber&=&\mathbb{E} \Big(\exp
i\sum_{j=1}^{m}\theta_{j}r^{-h}\int
(f(u+rt_{j},x)-f(u,x))d M(x)\Big)\\
\nonumber&=&\exp\Big\{-\int\big|\sum_{j=1}^{m}\theta_{j}r^{-h}(f(u+rt_{j},x)-f(u,x))\big|^{\alpha(x)}dx\Big\}\\
\nonumber
&=&\exp\Big\{-\int\big|\sum_{j=1}^{m}\theta_{j}r^{-h+1/\alpha(rz+u)}(f(u+rt_{j},rz+u)-f(u,rz+u))\big|^{\alpha(rz+u)}dz\Big\},\\
\label{thm-main-proof-1}
\end{eqnarray}
after setting $x=rz+u$.

Defining
\begin{equation*}\label{thm-main-proof-2}
Z(t)=\int h(t,z)dM_{\alpha(u)}(z),
\end{equation*}
its finite-dimensional distributions are given by the characteristic
function
\begin{equation}\label{thm-main-proof-3}
\mathbb{E}\Big(\exp i\sum_{j=1}^{m} \theta_{j}
Z(t_{j})\Big)=\exp\Big\{-\int\big|\sum_{j=1}^{m}\theta_{j}h(t_{j},z)\big|^{\alpha(u)}dz\Big\}.
\end{equation}
We now use Corollary $\ref{prop1}$, taking
\begin{equation*}\label{thm-main-proof-4}
g(r,z)=\sum_{j=1}^{m}\theta_{j}\frac{f(u+rt_{j},rz+u)-f(u,rz+u)}{r^{h-1/\alpha(rz+u)}},
\end{equation*}
\begin{equation*}\label{thm-main-proof-5}
k(z)=\sum_{j=1}^{m}\theta_{j}h(t_j,z),
\end{equation*}
and
\begin{equation*}\label{thm-main-proof-6}
\beta(x)=\alpha(u+x).
\end{equation*}
Then
\begin{equation*}\label{thm-main-proof-7}
\int|g(r,z)-k(z)|^{a,b}dz\\
\to 0,
\end{equation*}
as $r\to 0$, using  (\ref{thm-main-condition-1})and the quasi norm properties of $\|\cdot\|_a$ and $\|\cdot\|_b$. Thus by Corollary $\ref{prop1}$
\begin{equation*}\label{thm-main-proof-8}
\int\big|\sum_{j=1}^{m}\theta_{j}r^{-h+1/\alpha(rz+u)}(f(u+rt_{j},rz+u)-f(u,rz+u))\big|^{\alpha(rz+u)}dz\\
\to
\int\big|\sum_{j=1}^{m}\theta_{j}h(t_{j},z)\big|^{\alpha(u)}dz,
\end{equation*}
as $r\to 0$.

Since the exponential function is continuous,
($\ref{thm-main-proof-1}$), and hence (\ref{thm-main-proof-0}), is
 convergent to ($\ref{thm-main-proof-3}$) as $r\to 0$ for all $(\theta_1,\dots,\theta_m)$.
By L\'{e}vy's  Continuity Theorem, $r^{-h}(Y(u+rt)-Y(u))\fdd Z(t)$
as $r\to 0$, noting that (\ref{thm-main-proof-3}) is a characteristic function. Thus $Y$ is $h$-localisable with local form $Y_{u}'$ given by (\ref{thm-main-Yu}).

Finally, if (\ref{stongcond}) holds then by Proposition \ref{moments}, for $0<p<a$,
\begin{eqnarray*}
\mathbb{E}(|Y_r(t) - Y_r(v)|^p)
&=& \mathbb{E}\bigg(\Big|\int
\frac{f(u+rt,x)-f(u-rv,x)}{r^{h}}dM_\alpha(x)\Big|^p\bigg)\\
&\leq& 
c_2 \Big\|
\frac{f(u+rt,\cdot)-f(u-rv,\cdot)}{r^{h}}\Big\|_{\alpha}^p\\
&\leq& c_3 |t-v|^\eta p.
\end{eqnarray*}
By choosing $p$ such that $1/\eta<p< a$, Kolmogorov's continuity theorem, see \cite[Theorem 25.2]{RW}, implies that, for each $0 <\beta< (\eta a-1)/a$ and each bounded interval $I$, the process $Y_r$ satisfies an a.s. 
H\"{o}lder condition 
$$|Y_r(t) - Y_r(v)|\leq C_r |t-v|^\beta \quad (t,v \in I),$$
where the random constants behave uniformly in $r$, i.e,
 $\sup_{0<r \leq r_0}\mathbb{P} (C_r \geq m) \to 0$ as $m \to \infty$. Thus for all $\epsilon, \tau>0$ there exists $\delta>0$ such that
 $$\limsup_{r \to 0} \mathbb{P}\Big( \sup_{|t-v|<\delta,\; \;t,v \in I}
 |Y_r(t) - Y_r(v)|>\tau\Big) <\epsilon.$$ 
 In other words, the $Y_r$ are strongly stochastically equicontinuous on $I$ which, along with convergence of the finite dimensional distributions, implies that $Y_r$ converges to $Y'$ in distribution on the space of continuous functions with the metric of convergence on bounded intervals, see \cite[Theorem 8.2]{Bil} or\cite[Theorem 10.2]{Pol}  .
 $\Box$

\section{Examples}
\setcounter{equation}{0}
\setcounter{theo}{0}
We give  a number of examples to illustrate Theorem \ref{thm-main}. Some of these are considered in \cite{FL,FLL,LL} using alternative definitions of multistable processes.

It is convenient to make the convention that 
\begin{equation*}\label{indicator}
\mathbf{1}_{[u,v]}=-\mathbf{1}_{[v,u]},
\end{equation*}
if $v<u$ in the following examples.

\begin{exa}
Weighted multistable L\'{e}vy motion.

Let
\begin{equation*}\label{exa2-2}
Y(t)=\int w(x)\mathbf{1}_{[0,t]}(x)dM_{\alpha(x)}(x),
\end{equation*}
where $\alpha: \mathbb{R}\to [a,2] $ is continuous and $a>0$,
and $w: \mathbb{R} \to \mathbb{R}$ is continuous. Let $u \in \mathbb{R}$ be such that  $w(u) \neq 0$ and suppose
that as $v\to u$,
\begin{equation}
|\alpha(u)-\alpha(v)|=o\big(1/\big|\log|u-v|\big|\big). \label{logcond}
\end{equation}
Then $Y$ is $1/\alpha(u)$-localisable at $u$  with local form
\begin{equation*}\label{exa_makv_1}
Y_{u}'=\Big\{\int
w(u)\mathbf{1}_{[0,t]}(z)dM_{\alpha(u)}(z),\textrm{ }\textrm{
}\textrm{ }t\in\mathbb{R}\Big\}=w(u)L_{\alpha(u)},
\end{equation*}
where $L_{\alpha(u)}$ is a $\alpha(u)$-stable L\'{e}vy motion
\end{exa}

\begin{proof}
Take $f(t,x)=w(x)\mathbf{1}_{[0,t]}(x)$ and  $h(t,z)= w(u)\mathbf{1}_{[0,t]}(z)$. Condition (\ref{logcond}) ensures that 
$r^{1/\alpha(u)-1/\alpha(u+rz)}\to1$ as $r\to 0$ uniformly for $z\in [0,t]$ which is needed to ensure that (\ref{thm-main-condition-1}) holds. Then Theorem  \ref{thm-main} gives the conclusion.
\end{proof}

Next we consider multistable reverse Ornstein-Uhlenbeck motion. Notice that in the multistable case, we get a curious restriction on the range of $\alpha$.

\begin{exa}
Multistable reverse Ornstein-Uhlenbeck motion.

Let
\begin{equation}
Y(t)=\int_{t}^{\infty} \exp(-\lambda(x-t))dM_{\alpha(x)}(x),
\end{equation}
where $\alpha$: $\mathbb{R}\to [a,b] \subseteq (1,2]$ is continuous with $1<
\sqrt{b}< a \leq b\leq2$. Let $u\in \mathbb{R}$  and
suppose that as $v\to u$,
\begin{equation}
|\alpha(u)-\alpha(v)|=o\big(1/\big|\log|u-v|\big|\big).
\end{equation}
Then $Y$ is $1/\alpha(u)$-localisable at $u$ with local form
\begin{equation}\label{rou_1}
Y'_{u}=\Big\{\int -\mathbf{1}_{(0,t)}(z)dM_{\alpha(u)}(z),\textrm{
}\textrm{ }\textrm{ }t\in\mathbb{R}\Big\}.
\end{equation}
\end{exa}
\begin{proof}
We take  $f(t,x)=\exp(-\lambda(x-t))\mathbf{1}_{[t,\infty)}(x)$ and  $h(t,z)= -\mathbf{1}_{[0,t)}(z)$ in Theorem \ref{thm-main}. After a little simplification,
\begin{eqnarray*}
&&\hspace{-1cm}\int\Big|\frac{f(u+rt,u+rz)-f(u,u+rz)}{r^{1/\alpha(u)-1/\alpha(u+rz)}}-h(t,z)\Big|^{a,b}dz\\
&=&\int_{-|t|}^{|t|}\Big|\frac{-\exp(-\lambda
rz)\mathbf{1}_{[0,t)}(z)}{r^{1/\alpha(u)-1/\alpha(u+rz)}}+\mathbf{1}_{[0,t)}(z)\Big|^{a,b}dz
+\int_{|t|}^{\infty}\Big|\frac{\exp(-\lambda rz)(\exp(\lambda
rt)-1)\mathbf{1}_{[t,\infty)}(z)}{r^{1/\alpha(u)-1/\alpha(u+rz)}}\Big|^{a,b}dz.
\end{eqnarray*}
The first integral converges to $0$, noting that 
$r^{1/\alpha(u)-1/\alpha(u+rz)}\to 1$ as $r\to 0$, uniformly on $z \in [-t,t]$. The second integral is bounded by
\begin{eqnarray*}
&&\hspace{-2cm}\int_{|t|}^{\infty}\big|r^{-1/a+1/b}\exp(-\lambda
rz)(\exp(\lambda rt)-1)\big|^{a,b}dz\\
&\leq& r^{1-b/a}\int_{|t|}^{\infty}\big|\exp(-\lambda
rz)(\exp(\lambda rt)-1)\big|^{a,b}dz\\
&\leq&c_{1}r^{1-b/a}|\exp(\lambda
rt)-1|^{a}\int_{|t|}^{\infty}\big|\exp(-\lambda arz)\big|dz\\
&\leq& c_{2}r^{1-b/a}(\lambda r|t|)^{a}\exp (-\lambda ar|t|)(\lambda
ra)^{-1}\\
&\leq& c_{3} r^{a-b/a},
\end{eqnarray*}
for fixed $t$, where $c_{1}$, $c_{2}$ and $c_{3}$ are independent of $r<1$.
Since $a-b/a >0$ the second integral converges to $0$, so the conclusion follows from 
 Theorem \ref{thm-main}.
\end{proof}

The next example is linear fractional multistable motion. Recall from \cite{Bk_Sam} that asymmetric linear fractional $\alpha$-stable motion is given by
\begin{equation}\label{defi-linfracalphastab-1}
L_{\alpha,h,b^{+},b^{-}}(t)=\int_{-\infty}^{\infty}\varrho_{\alpha,h}(b^{+},b^{-},t,x)dM_{\alpha}(x)
\end{equation}
where $t,b^{+},b^{-}\in \mathbb{R}$, and
\begin{equation*}
\varrho_{\alpha,h}(b^{+},b^{-},t,x)=b^{+}\big((t-x)_{+}^{h-1/\alpha}-(-x)_{+}^{h-1/\alpha}\big)+b^{-}\big((t-x)_{-}^{h-1/\alpha}-(-x)_{-}^{h-1/\alpha}\big),
\end{equation*}
and $M_{\alpha}$ is  $\alpha$-stable random measure
$(0<\alpha<2)$. By convention, if $h-1/\alpha =0$, we take 
$$\varrho_{\alpha,h}(b^{+},b^{-},t,x) =(b^+ -b^-)\mathbf{1}_{[0,t]}(x)$$
if $t\geq 0$, and
$$\varrho_{\alpha,h}(b^{+},b^{-},t,x) =-(b^+ -b^-)\mathbf{1}_{[t,0]}(x)$$
if $t <0$. Then (\ref{defi-linfracalphastab-1}) is an $\alpha$-stable process. 

For a multistable version let $\alpha: \mathbb{R} \to [a,b] \subseteq (0,2)$ be continuous. We define {\em linear fractional $\alpha(x)$-multistable motion} by
\begin{equation}\label{defi-linfracalphastab-2}
L_{\alpha(x),h,b^{+},b^{-}}(t)=\int_{-\infty}^{\infty}\varrho_{\alpha(x),h}(b^{+},b^{-},t,x)dM_{\alpha(x)}(x)
\end{equation}
where $t\in\mathbb{R}$, $b^{+},b^{-}\in \mathbb{R}$, and
\begin{equation*}
\varrho_{\alpha(x),h}(b^{+},b^{-},t,x)=b^{+}\big((t-x)_{+}^{h-1/\alpha(x)}-(-x)_{+}^{h-1/\alpha(x)}\big)+b^{-}\big((t-x)_{-}^{h-1/\alpha(x)}-(-x)_{-}^{h-1/\alpha(x)}\big),
\end{equation*}
where  $M_{\alpha}(x)$ is $\alpha(x)$-multistable random measure.

It may be checked directly that if $t \in \mathbb{R}$ and  $1/a-1/b<h<1+1/b-1/a$ then 
$\varrho_{\alpha(\cdot),h}(b^{+},b^{-},t,.) \in \mathcal{F}_{a,b}$ so that (\ref{defi-linfracalphastab-2}) is well-defined.

We show that linear fractional multistable motion has linear stable motion as its local form. We consider the case when $b^{+}=1$ and $b^{-}=0$, the argument is similar for other $b^{+}$ and $b^{-}$.

\begin{prop}\label{Ch3_exa_lin}
Linear fractional multistable motion.

Let
\begin{eqnarray*}
\nonumber Y(t)&=&\int
(t-x)_{+}^{h-1/\alpha(x)}-(-x)_{+}^{h-1/\alpha(x)}dM_{\alpha(x)}(x) \\
\nonumber &=&\int \varrho_{\alpha(x),h}(1,0,t,x)dM_{\alpha(x)}(x)\\
&=&L_{\alpha(x),h,1,0}(t),
\end{eqnarray*}
where $\alpha$: $\mathbb{R}\to [a,b] \subseteq (0,2)$ is continuous.
If 
\begin{equation}\label{exa_lin_12} 
1/a-1/b<h<1+1/b-1/a, 
\end{equation}
then $Y$ is $h$-localisable at each $u \in \mathbb{R}$ with
local form
\begin{eqnarray*}
\nonumber
Y_{u}'(t)&=&\Big\{\int\big((t-z)_{+}^{h-1/\alpha(u)}-(-z)_{ +}^{h-1/\alpha(u)}\big)dM_{\alpha(u)}(z),\textrm{
}\textrm{ }\textrm{ }t\in \mathbb{R}\Big\}\\
&=& L_{\alpha(u),h,1,0}(t).
\end{eqnarray*}
Furthermore, if $1/a<h<1+1/b-1/a$, then $Y$ has a continuous version and is strongly localisable at each $u \in \mathbb{R}$.
\end{prop}

\begin{proof}
We take
$f(t,x)=(t-x)_{+}^{h-1/\alpha(x)}-(-x)_{+}^{h-1/\alpha(x)} \in \mathcal{F}_{a,b}$, given (\ref{exa_lin_12}),
and $h(t,z)=(t-z)_{+}^{h-1/\alpha(u)}-(-z)_{+}^{h-1/\alpha(u)}$ in Theorem \ref{thm-main}.   Then\begin{eqnarray}
\nonumber &&\hspace{-1cm}\int\Big|\frac{f(u+rt,u+rz)-f(u,u+rz)}{r^{h-1/\alpha(u+rz)}}-h(t,z)\Big|^{a,b}dz\\
\nonumber
&=&\int\Big|(t-z)_{+}^{h-1/\alpha(u+rz)}-(-z)_{+}^{h-1/\alpha(u+rz)}-(t-z)_{+}^{h-1/\alpha(u)}+(-z)_{+}^{h-1/\alpha(u)}\Big|^{a,b}dz.
\end{eqnarray}
This integral converges to $0$ as $r \to 0$. This may be established by breaking the range of integration in the parts: $|z| < \delta, |z-t| < \delta, |z|>M$ and
$A=\{z: \delta \leq |z| \leq M$ and $\delta \leq |z-t|\}$.
By choosing sufficiently small $\delta$ and large $M$, the integral over the first three parts can be made arbitrarily small, uniformly as $r \to 0$.  The integrand converges to $0$ pointwise on $A$ and the bounded convergence theorem gives the integral over $A$ convergent to $0$. The conclusion follows from Theorem \ref{thm-main}. 

Finally, if $1/a < h <1$ it is easily checked by a routine integral estimate (if $t>u$ splitting the resulting integral at $u$) that (\ref{ctsverequiv}) holds so $Y$ has a continuous version, and similarly that (\ref{equivstongcond}) if $1/a < \eta <h$, so $Y$ is strongly localisable by Theorem \ref{thm-main}. 
\end{proof}

\bibliographystyle{plain}

\end{document}